\documentclass[12pt]{article}
\usepackage{amsfonts}
\usepackage{latexsym,amsmath,amssymb}
\usepackage{graphics,epstopdf}
\usepackage[pdftex]{graphicx}
\newtheorem{definition}{Definition}[section]
\newtheorem{theorem}[definition]{Theorem}

\newtheorem{corollary}[definition]{Corollary}

\numberwithin{equation}{section}

\begin{document}

\begin{center}
{\Large \textbf{\ Approximation by $(p,q)$-Lorentz polynomials on a compact disk}}

\bigskip

\textbf{M. Mursaleen}, \textbf{Faisal Khan}, \textbf{Asif Khan}

\bigskip

\textbf{Abstract}
\end{center}

\parindent=8mm {\footnotesize {In this paper, we introduce a new analogue of Lorentz
polynomials based on $(p,q)$-integers and we call it as $(p,q)$-Lorentz polynomials. We obtain quantitative
estimate in the Voronovskaja's type thoerem  and exact orders in simultaneous approximation
by the complex $(p,q)$-Lorentz polynomials of degree $n\in \mathbb{N},$ where $q>p>1$ attached
to analytic functions in compact disks of the complex plane.}}

\bigskip

\parindent=8mm{\footnotesize \emph{AMS Subject Classifications (2010)}: {41A10, 41A25, 41A36.}}

\parindent=8mm{\footnotesize \emph{Keywords and phrases}: $(p,q)$-integer; Complex $(p,q)$-Lorentz polynomials;
Voronovskaja's theorem; iterates; compact disk.}

\section{ Introduction and preliminaries}

\hspace{8mm}In 1986, G. G. Lorentz \cite{lgga}, introduced the following sequence of operators
defined for any analytic function $f$ in a domain containing the origin
\begin{equation}\label{e1.1}
L_{n}(f;z)=\sum\limits_{k=0}^{n}\left(
\begin{array}{c}
n \\
k%
\end{array}%
\right)
\left(\frac{z}{n}
\right)^k
f^{(k)}(0),~~n\in \mathbb{N}.
\end{equation}

\parindent=8mm In last two decades, the applications of $q$-calculus emerged as a new area
in the field of approximation theory. 

In \cite{sggo}, Gal introduced and studied the $q$-analouge of the Lorentz operators for $q>1,$ for
any analytic function $f$ in a domain containing the origin as follows:

\begin{equation}\label{e1.2}
L_{n,q}(f;z)=\sum\limits_{k=0}^{n}q^{\frac{k(k-1)}{2}}\left[
\begin{array}{c}
n \\
k%
\end{array}%
\right]_q
\left(\frac{z}{[n]_q}
\right)^k
D^{(k)}_q(f)(0),~~n\in \mathbb{N},~z\in\mathbb{C}
\end{equation}
Details on the $q$-calculus can be found in \cite{gp1,vp}

Several authors have introduced and studied the approximation properties for different operators
in compact disk. For instance, in \cite{m1,m2,mk} Mahmudov studied $q$-Stancu polynomials,
q-Szasz Mirakjan operators and generalised Kantorovich operators; In
\cite{gmk} Gal et al studied q-szasz-Kantorovich operators.

Recently, Mursaleen et al applied $(p,q)$-calculus in approximation theory and introduced first
$(p,q)$-analogue of Bernstein operators \cite{mka1}. Similarly they introduced and studied
approximation properties for $(p,q)$-Bernstein-Stancu operators {\cite{mka2}} and
$(p,q)$-Bernstein-Kantorovich operators {\cite{mka3}}.

Let us recall certain notations of $(p,q)$-calculus.

The $(p,q)$ integer $[n]_{p,q}$ is defined by
\begin{equation*}
[n]_{p,q}:=\frac{p^n-q^n}{p-q},~~~n=0,1,2,\cdots,~~p>q\geq1.
\end{equation*}
The $(p,q)$-binomial expansion is given as
\begin{equation*}
(ax+by)_{p,q}^{n}:=\sum\limits_{k=0}^{n}\left[
\begin{array}{c}
n \\
k%
\end{array}%
\right] _{p,q}a^{n-k}b^{k}x^{n-k}y^{k}
\end{equation*}%
\begin{equation*}
(x+y)_{p,q}^n:=(x+y)(px+qy)(p^2x+q^2y)\cdots(p^{n-1}x+q^{n-1}y),
\end{equation*}%
and the $(p,q)$-binomial coefficients are defined by
\begin{equation*}
\left[
\begin{array}{c}
n \\
k%
\end{array}%
\right] _{p,q}:=\frac{[n]_{p,q}!}{[k]_{p,q}![n-k]_{p,q}!}.
\end{equation*}

The $(p,q)$-derivative of the function $f$ is defined as
\begin{equation*}
D_{p,q}f(x)=\frac{f(px)-f(qx)}{(p-q)x},~~~x\neq0,
\end{equation*}

\noindent
and $(D_{p,q})f(0)=f^{\prime}(0),$ provided that $f$ is differentiable at $0.$
It can be easily\\
seen that $D_{p,q}x^n=[n]_{p,q}x^{n-1}.$\\
\noindent
Details on $(p,q)$-calculus can be found in \cite{jacob,jag,mah,sad,vivek}.

 Our results generalize the results of Gal \cite{sggo}

\section{Construction of Operators}

\hspace{8mm}Now, with the  help of $(p,q)$-calculus and using above formula, we present $(p,q)$-analogue of
Lorentz operators  (\ref{e1.1}) as follows:
\begin{equation}\label{e2.2}
L_{n,p,q}(f;z)=\sum\limits_{k=0}^{n}q^{\frac{k(k-1)}{2}}\left[
\begin{array}{c}
n \\
k%
\end{array}%
\right]_{p,q}
\left(\frac{z}{[n]_{p,q}}
\right)^k
D^{(k)}_{p,q}(f)(0),~~n\in \mathbb{N},~z\in\mathbb{C}%
\end{equation}
Note that for $p=1,$ $(p,q)$-Lorentz operators given by (\ref{e2.2}) turn out to
be $q$-analogue of  Lorentz operators (\ref{e1.2}).

Firstly, we obtain an upper approximation estimate.

\begin{theorem}\label{theorem2.1} {\em Let $R>q>p>1$ and $D_{R}=\{z\in\mathbb{C}: |z|<R\}.$ Suppose that
$f:D_R\rightarrow\mathbb{C}$ is analytic in $D_R$, i.e., $f(z)=\sum\limits_{k=0}^{\infty}c_kz^k,$ for all $z\in D_R.$
 \begin{enumerate}
 \item [(i)] Let $1\leq r < \frac{pr_1}{q}<\frac{pR}{q}$ be arbitrary fixed. For all $|z|\leq r$ and $n\in\mathbb{N}.$
Then, we have the upper estimate as
 \begin{equation*}
\left|L_{n,p,q}(f)(z)-f(z)\right|\leq \frac{p^n}{[n]_{p,q}}~M_{r_1,p,q}(f),
\end{equation*}
where $M_{r_1,p,q}(f)=\frac{p(q-p+1)}{(q-p)^2}\sum\limits_{k=0}^{\infty}|c_k|(k+1)r^k_1<\infty.$
\item [(ii)] Let $1\leq r < <r^{\ast}<\frac{pr_1}{q}<\frac{pR}{q}$ be arbitrary fixed. For the simultaneous
approximation by complex Lorentz polynomials, for all $|z|\leq r,~~m,n\in\mathbb{N},$ we have
 \begin{equation*}
\left|L_{n,p,q}^{(m)}(f)(z)-f^{(m)}(z)\right|\leq \frac{p^n}{[n]_{p,q}}~M_{r_1,p,q}(f)\frac{m!~r^{\ast}}{(r^{\ast}-r)^{m+1}},
\end{equation*}
where $M_{r_1,p,q}(f)$ is given as at the above point (i).
\end{enumerate}}
\end{theorem}
\noindent\textbf{Proof.} (i) For $e_j(z)=z^j,$ it is to see that
$L_{n,p,q}(e_0)(z)=1,~L_{n,p,q}(e_1)(z)=z.$ Then we have
\begin{equation*}
 L_{n,p,q}(e_j)(z)=q^{\frac{j(j-1)}{2}}
\left[
\begin{array}{c}
n \\
j
\end{array}%
\right] _{p,q}
[j]_{p,q}!~\frac{z^j}{[n]_{p,q}^j}, ~~2 \leq j \leq n, ~\mbox{for all}~ j,n\in\mathbb{N}
\end{equation*}
by some simple calculation, we get 
 \begin{equation*}
 L_{n,p,q}(e_j)(z)=z^j~\left(1-p^{n-1}\frac{[1]_{p,q}}{[n]_{p,q}}\right)
\left(1-p^{n-2}\frac{[2]_{p,q}}{[n]_{p,q}}\right)...\left(1-p^{n-(j-1)}\frac{[j-1]_{p,q}}{[n]_{p,q}}\right).
\end{equation*}
It is easy to see that for $j\geq n+1,$ we get $ L_{n,p,q}(e_j)(z)=0.$\\
Now it can be easily seen that
 \begin{equation*}
 L_{n,p,q}(f)(z)=\sum\limits_{j=0}^{\infty}c_j L_{n,p,q}(e_j)(z)~~\mbox{for all}~~|z|\leq r.
\end{equation*}
Hence

\vspace{2mm}
\noindent
$\big| L_{n,p,q}(f)(z)-f(z)\big|$
\begin{align*}
{}~~~~~ &\leq\sum\limits_{j=0}^n|c_j|\big| L_{n,p,q}(e_j)(z)-e_j(z)\big|+\sum\limits_{j=n+1}^{\infty} |c_j|
 \big|L_{n,p,q}(e_j)(z)-e_j(z)\big|\\
 &\leq\sum\limits_{j=2}^n|c_j|~r^j~\frac{1}{p^{jn}}~\left|\left(1-p^{n-1}\frac{[1]_{p,q}}{[n]_{p,q}}\right)
\left(1-p^{n-2}\frac{[2]_{p,q}}{[n]_{p,q}}\right)...\left(1-p^{n-(j-1)}\frac{[j-1]_{p,q}}{[n]_{p,q}}\right)-1\right|\\
 & ~~~+\sum\limits_{j=n+1}^{\infty} |c_j| r^j,~~~\mbox{for all}~~|z|\leq r.
\end{align*}
%

On the other hand, the analyticity of $f$ implies $c_j=\frac{f^{(k)}{(0)}}{j!},$ and by the Cauchy's estimates
of the coefficients $c_j$ in the disk $|z|\leq r_1,$ we have $|c_j|\leq\frac{K_{r_1}}{r-1},$ for all $j\geq 0,$
where
\begin{equation*}
K_{r_1}=\max\{|f(x)|:|z|\leq r_1\}\leq\sum\limits_{j=2}^n|c_j|r^j\leq\sum\limits_{j=2}^n|c_j|~(j+1)~r_1^j:=R_{r_1}(f)<\infty.
\end{equation*}
Therefore
\begin{equation*}
\sum\limits_{j=n+1}^{\infty}|c_j|r^j\leq R_{r_1}(f)\left[\frac{r}{r_1}\right]^{n+1}
\sum\limits_{j=0}^{\infty}\left(\frac{r}{r_1}\right)^j=R_{r_1}(f)\left[\frac{r}{r_1}\right]^{n+1}.\frac{r_1}{r_1-r}
\end{equation*}
\begin{equation*}
  =R_{r_1}(f)\frac{r}{r_1-r}\left[\frac{r}{r_1}\right]^n\leq R_{r_1}(f)\frac{p^{n+1}}{[n]_{p,q}}\frac{1}{(q-p)^2},
\end{equation*}
and finally we get
\begin{equation*}
 \big| L_{n,p,q}(f)(z)-f(z)\big|\leq\frac{p^{n+1}}{[n]_{p,q}}\frac{(q-p+1)}{(q-p)^2}R_{r_1}(f)
\end{equation*}
for all $n\in \mathbb{N}$ and $|z|\leq r.$\\

\noindent
(ii) Let $\gamma$ be the circle of radius $r^{\ast}>r$ and center $0,$ since
for any $|z|\leq r$ and $\upsilon\in \gamma,$ we have $|\upsilon-z|\geq r^{\ast}-r.$
By Cauchy's formula it follows that for all $n\in\mathbb{N}$
\begin{align*}
\big| L_{n,p,q}^{(m)}(f)(z)-f^{(m)}(z)\big|& =\frac{m!}{2\pi}\left|\int_{\gamma} \frac{L_{n,p,q}(f)(\upsilon)-f(\upsilon)}{(\upsilon-z)^{m+1}}d\upsilon\right|\\
&\leq\frac{{p^{n+1}}}{{[n]_{p,q}}}~M_{r_1,p,q}(f)~\frac{m!}{2\pi}~\frac{2\pi r^{\ast}}{(r^{\ast}-r)^{m+1}}\\
&=\frac{p^{n+1}}{[n]_{p,q}}~M_{r_1,p,q}(f)~\frac{m!~r^{\ast}}{(r^{\ast}-r)^{m+1}}.
\end{align*}
We have the following quantitative Voronovskaja-type results.

\begin{theorem}\label{theorem2.2} {\em For $R>q^4>p^4>1,$ let $f:D_R\rightarrow\mathbb{C}$
be analytic in $D_R$, i.e., $f(z)=\sum\limits_{k=0}^{\infty}c_kz^k,$ for all $z\in D_R$ and
let $1\leq r<\frac{p^3r_1}{q^3}<\frac{p^4R}{q^4}$ be arbitrary fixed. Then for all
$n\in\mathbb{N},~|z|\leq r,$ we have
\begin{equation*}
\left|L_{n,p,q}(f)(z)-f(z)+\frac{S_{p,q}(f)(z)}{[n]_{p,q}}\right|\leq\frac{p^{2n}}{{[n]_{p,q}^2}}~Q_{r_1,p,q}(f),
\end{equation*}
where
\begin{equation*}
S_{p,q}(f)(z)=\sum\limits_{k=2}^{\infty}p^{n-(k-1)}c_k\frac{[k]_{p,q}-[k]_q}{q-1}z^k
=\sum\limits_{k=2}^{\infty}p^{n-(k-1)}c_k\big([1]_{p,q}+...+[k-1]_{p,q}\big)
\end{equation*}
and $Q_{r_1,p,q}=\frac{pq-q+p}{(p-1)(q-p)}~\sum\limits|c_k|(k+1)(k+2)^2\left(\frac{q}{p}r_1\right)^k<\infty.$}
\end{theorem}
\noindent\textbf{Proof.} We have
\begin{equation*}
\left |L_{n,p,q}(f)(z)-f(z)+\frac{S_{p,q}(f)(z)}{[n]_{p,q}}\right|\hspace{8cm}
\end{equation*}
\begin{align*}
{}&=\left|\sum\limits_{k=0}^{\infty}c_k\left [L_{n,p,q}(e_k)(z)-e_k(z)+p^{n-(k-1)}
~\frac{[k]_{p,q}-[k]_q}{p-1}~e_k(z)\right]\right| \\
&\leq\left|\sum\limits_{k=0}^{n}c_k\left [L_{n,p,q}(e_k)(z)-e_k(z)+p^{n-(k-1)}
~\frac{[k]_{p,q}-[k]_q}{p-1}~e_k(z)\right]\right|\\
&\hspace{2cm}+\left|\sum\limits_{k=n+1}^{\infty}p^{n-(k-1)}~c_k~ z^k\left(\frac{[k]_{p,q}-[k]_q}{p-1}-1\right)\right|\\
&\leq\left|\sum\limits_{k=0}^{n}c_k\left [L_{n,p,q}(e_k)(z)-e_k(z)+p^{n-(k-1)}
~\frac{[k]_{p,q}-[k]_q}{p-1}~e_k(z)\right]\right|\\
&\hspace{2cm}+\sum\limits_{k=n+1}^{\infty}~|c_k|~ r^k\left(p^{n-(k-1)}~\frac{[k]_{p,q}-[k]_q}{p-1}-1\right),\\
\end{align*}
for all $|z|\leq r$ and $n\in\mathbb{N}.$\\
In what follows, firstly we will prove by mathematical induction with respect to $k$ that
\begin{equation}\label{e2.3}
0 \leq E_{n,k,p,q}(z) \leq \frac{p^{2n}}{[n]_{p,q}^2}\frac{(k+1)(k-2)^2}{(q-p)}\left(\frac{qr_1}{p}\right)^{k},
\end{equation}
for all $2\leq k\leq n$ (here $n\in\mathbb{N}$ is arbitrary fixed) and $|z|\leq r,$ where
\begin{equation*}
E_{n,k,p,q}(z)=L_{n,p,q}(e_k)(z)-e_k(z)+\frac{p^{n-(k-1)}}{[n]_{p,q}}~\frac{[k]_{p,q}-[k]_q}{p-1}~e_k(z)
\end{equation*}
\begin{equation*}
=L_{n,p,q}(e_k)(z)-e_k(z)+\frac{p^{n-(k-1)}}{[n]_{p,q}}\left([1]_{p,q}+...+[k-1]_{p,q}\right)~e_k(z).
\end{equation*}
By mathematical induction, we easily
\begin{equation*}
  \frac{[k]_{p,q}-[k]_q}{p-1}=\left([1]_{p,q}+...+[k-1]_{p,q}\right).
\end{equation*}

On the other hand, by the formula for $L_{n,p,q}(e_k)$ in the proof of Theorem \ref{theorem2.1} (i),
simple calculation leads to $E_{n,2,p,q}(z)=0,$ for all $n\in\mathbb{N}$ and to
the recurrence relation
\begin{equation*}
  E_{n,k,p,q}(z)=-~\frac{z^2}{[n]_{p,q}}p^{n-(k-1)}~D_{p,q}[L_{n,p,q}(e_{k-1})(z)-e_{k-1}(z)]
  ~+~\frac{p-1}{p}z[L_{n,p,q}(e_{k-1})(z)-e_{k-1}(z)]
\end{equation*}
\begin{equation*}
  \hspace{10cm}+~\frac{z}{p}~E_{n,{k-1},p,q}(z),~~|z|\leq r.
\end{equation*}
Now, for $|z|\leq r$ and $3\leq k \leq n$ and applying the
mean value theorem in complex analysis, with  notation $\|f\|_r=\max\{|f(z)|:|z|\leq r\},$ we get
\begin{equation*}
  \left|E_{n,k,p,q}(z)\right|=\frac{r^2}{[n]_{p,q}}~p^{n-(k-1)}
  ~\|\big(L_{n,p,q}(e_{k-1})(z)-e_{k-1}(z)\big)^{\prime}\|_{\frac{qr}{p}}\hspace{4cm}
\end{equation*}

\begin{equation*}
 |E_{n,k,p,q}(z)|\leq \left(\frac{p^n}{[n]_{p,q}}(k+1)\right)~~\frac{p^n}{[n]_{p,q}}(k-2)[k-2]_{p,q}~r^k_1+~r_1\left|E_{n,{k-1},p,q}(z)\right|
\end{equation*}
\begin{equation*}
 \leq \frac{p^{2n}}{[n]^2_{p,q}}(k+1)(k-2)[k-2]_{p,q}r_1^k~+~r_1\left|E_{n,{k-1},p,q}(z)\right|
\end{equation*}
Now on taking $k=1,2,3,...,$ step by step, we easily obtain the estimate
\begin{equation*}
\left|E_{n,k,p,q}(z)\right|\leq \frac{p^{2n}}{[n]^2_{p,q}}r_1^k\sum\limits_{j=3}^k(j-1)(j-2)[j-2]_{p,q}
\end{equation*}
\begin{equation*}
  \leq \frac{p^{2n}}{[n]_{p,q}^2}\frac{(k+1)(k-2)^2}{(q-p)}\left(\frac{qr_1}{p}\right)^{k}
\end{equation*}
Now we calculate\\

$\left|\sum\limits_{k=0}^{n}|c_k|\left [L_{n,p,q}(e_k)(z)-e_k(z)+q^{n-(k-1)}~\frac{[k]_{p,q}-[k]_q}{p-1}~e_k(z)\right]\right|\hspace{2cm}$
\begin{align*}
{}&\hspace{1cm}\leq\sum\limits_{k=0}^n~|c_k|~\left|E_{n,k,p,q}(z)\right| \\
&\leq\frac{p^{2n}}{[n]_{p,q}^2}\frac{1}{(q-p)}\sum\limits_{k=0}^{n}|c_k|(k+1)(k-2)^2\left(\frac{qr_1}{p}\right)^{k}\\
&\leq\frac{p^{2n}}{[n]_{p,q}^2}\frac{1}{(q-p)}\sum\limits_{k=0}^{n}|c_k|(k+1)(k+2)^2\left(\frac{qr_1}{p}\right)^{k}\\
\end{align*}
On the other hand, since $\left(p^{n-(k-1)}~\frac{[k]_{p,q}-[k]_q}{p-1}-1\right)\geq 0$ for all $k\geq n+1,$ similar
to proof of Theorem \ref{theorem2.1} (i), we get
\begin{equation*}
  \sum\limits_{k=n+1}^{\infty}~|c_k|~ r^k\left(p^{n-(k-1)}~\frac{[k]_{p,q}-[k]_q}{p-1}-1\right)\leq \sum\limits_{k=n+1}^{\infty}p^{n-(k-1)}~|c_k|~r^k\frac{[k]_{p,q}}{(p-1)[n]_{p,q}}
\end{equation*}
\begin{equation*}
  \leq\sum\limits_{k=n+1}^{\infty}p^{n-(k-1)}~|c_k|~\frac{1}{(p-1)[n]_{p,q}}\frac{k q^k}{p^k(q-p)}\hspace{0.5cm}
\end{equation*}
\begin{equation*}
  \leq\frac{R_{r_1}(f)p^{n+1}}{(p-1)[n]_{p,q}}\sum\limits_{k=n+1}^{\infty}~\frac{r^k}{r_1^k}~q^{k}p^{-k}\hspace{1.5cm}
\end{equation*}
\begin{equation*}
  \leq\frac{R_{r_1}(f)p^{n+1}}{(p-1)[n]_{p,q}}\sum\limits_{k=n+1}^{\infty}\left[\left(\frac{r}{r_1}\right)^{1/3}\right]^k
  ~\left[\left(\frac{r}{r_1}\right)^{1/3}\right]^{2k}q^{k}p^{-k}
\end{equation*}
\begin{equation*}
  \leq\frac{R_{r_1}(f)p^{n+1}}{(p-1)[n]_{p,q}}\left(\frac{r}{r_1}\right)^{\frac{(n+1)}{3}}
  \sum\limits_{k=0}^{\infty}~\left[\left(\frac{r}{r_1}\right)^{1/3}\right]^k
\end{equation*}
\begin{equation*}
  =\frac{R_{r_1}(f)p^{n+1}}{(p-1)[n]_{p,q}}\left(\frac{r}{r_1}\right)^{\frac{n}{3}}
  \frac{r^{\frac{1}{3}}}{(r_1^{\frac{1}{3}}-r^{\frac{1}{3}})}
\end{equation*}
\begin{equation*}
  \leq\frac{p^{2n+2}}{[n]^2_{p,q}}\frac{R_{r_1}(f)}{(p-1)~(q-p)^2}\hspace{3cm}
\end{equation*}
\begin{equation}\label{e2.3}
  \leq\frac{p^{2n+2}}{[n]^2_{p,q}~(p-1)~(q-p)^2}~\sum\limits_{k=0}^{n}|c_k|~(k+1)(k+2)^2~\left(\frac{q}{p}r_1\right)^k,
\end{equation}
where we used the inequalities, $[k]_{p,q}\leq\frac{kq^k}{p^k},~\frac{p^n}{q^n}\leq\frac{p^n}{(q-p)[n]_{p,q}}$ and $\frac{r^{1/3}}{(r_1^{1/3}-r^{1/3})}\leq \frac{p}{(q-p)}.$
Hence, by combining all above estimates, we have
\begin{equation*}
\left |L_{n,p,q}(f)(z)-f(z)+\frac{S_{p,q}(f)(z)}{[n]_{p,q}}\right|\
\leq\frac{(pq-q+p-1)}{(p-1)~(q-p)^2}~\frac{p^{2n}}{[n]^2_{p,q}}~\sum\limits_{k=0}^{n}|c_k|~(k+1)(k+2)^2~\left(\frac{q}{p}r_1\right)^k.
\end{equation*}
The following result gives the lower approximation estimate
\begin{theorem}\label{theorem2.3} {\em Let $R>p^4/q^4>1,~ f:D_R\rightarrow\mathbb{C}$
be analytic in $D_R$, i.e., $f(z)=\sum\limits_{k=0}^{\infty}c_kz^k,$ for all $z\in D_R$ and
let $1\leq r<\frac{q^3r_1}{p^3}<\frac{q^4R}{p^4}$ be arbitrary fixed. If $f$ is not a polynomial of degree $\leq 1,$ then for all
$n\in\mathbb{N}~\mbox{and}~|z|\leq r,$ we have
\begin{equation*}
  \left\| L_{n,p,q}(f)-f \right\|_{r}\geq\frac{p^n}{[n]_{p,q}}C_{r,r_1,p,q}(f),
\end{equation*}
 where the constant $C_{r,r_1,p,q}(f)$ depends only on $f,~r$ and $r_1.$
Here $\|f\|_r$ denotes $\max\limits_{|z|\leq r}\{|f(z)|\}.$}
\end{theorem}

\noindent\textbf{Proof.} For $S_{n,p,q}(f)(z)$ as defined in  Theorem \ref{theorem2.3}, all $|z|\leq r $
and $n\in\mathbb{N},$ we have

\vspace{8mm}
$L_{n,p,q}(f)(z)-f(z)$
\begin{equation*}
  \hspace{2cm}=\frac{p^n}{[n]_{p,q}}\left\{-S_{p,q}(f)(z)+\frac{p^n}{[n]_{p,q}}\left[\frac{[n]^2_{p,q}}{p^{2n}}
  \left(L_{n,p,q}(f)(z)-f(z)+\frac{p^n}{[n]_{p,q}}S_{p,q}(f)(z)\right)\right]\right\}.
\end{equation*}
Using the  following inequality
\begin{equation*}
  \|F+G\|_r\geq \big|\|F\|_r-\|G\|_r\big|\geq \|F\|_r-\|G\|_r.
\end{equation*}
We have\\
$ \|L_{n,p,q}(f)-f\|_r$
\begin{equation*}
  \hspace{2cm}\geq\frac{p^n}{[n]_{p,q}}\left\{\|S_{p,q}(f)(z)\|-\frac{p^n}{[n]_{p,q}}\left[\frac{[n]^2_{p,q}}{p^{2n}}
  \left\|L_{n,p,q}(f)(z)-f(z)+\frac{p^n}{[n]_{p,q}}S_{p,q}(f)(z)\right\|_r\right]\right\}.
\end{equation*}
Since by hypothesis $f$ is not a polynomial of degree $\leq 1$ in $D_R,$ we get $\|S_{p,q}(f)\|_{r}>0.$

Indeed, supposing the contrary it follows that $S_{p,q}(f)(z)=0$ for all $z\in\overline{{D_R}}=\{z\in \mathbb{C}:|z|\leq r\}.$

A simple calculation yields $S_{p,q}(f)(z)=z~\frac{D_{p,q}(f)(z)-f^{\prime}(z)}{p-q},~~S_{p,q}(f)(z)=0$
implies that $D_{p,q}(f)(z)=f^{\prime}(z),$ for all $z\in\overline{{D_r}}\backslash \{0\}.$ Taking into account
the representation of $f$ as $f(z)=\sum\limits_{k=0}^{\infty}c_kz^k,$ the last inequality immediately leads to $c_k=0,$
for all $k\geq 2,$ which means that $f$ is linear in $\overline{{D_r}},$ a contradiction with hypothesis.

Now, by Theorem \ref{theorem2.2} we have
\begin{equation*}
  \frac{[n]^2_{p,q}}{p^{2n}}
  \left\|L_{n,p,q}(f)(z)-f(z)+\frac{p^n}{[n]_{p,q}}S_{p,q}(f)(z)\right\|_r\leq Q_{r_1,p,q}(f),
\end{equation*}
where $ Q_{r_1,p,q}(f)$ is a positive constant depending only on $f,~r_1,~p$ and $q.$\\
Since $\frac{p^n}{[n_{p,q}]}\rightarrow0$ as $n\rightarrow \infty,$ there exists an index $n_{\circ}$
depending only on $f,~r,~r_1,~p$ and $q$ such that for all $n>n_{\circ},$ we have
\begin{equation*}
  \|S_{p,q}(f)(z)\|-\frac{p^n}{[n]_{p,q}}\left[\frac{[n]^2_{p,q}}{p^{2n}}
  \left\|L_{n,p,q}(f)(z)-f(z)+\frac{p^n}{[n]_{p,q}}S_{p,q}(f)(z)\right\|_r\right]\geq \frac{1}{2}\left\|S_{p,q}(f)\right\|_r,
\end{equation*}
which immediately implies that
\begin{equation*}
  \|L_{n,p,q}(f)-f\|_r\geq\frac{p^n}{[n]_{p,q}}\frac{1}{2} \|S_{p,q}(f)(z)\|_r,~~\mbox{for all}~n>n_{\circ}.
\end{equation*}

For $n\in\{1,...,n_{\circ}\},$ we have $\|L_{n,p,q}(f)-f\|_r~\geq~\frac{p^n}{[n]_{p,q}}~M_{r,r_1,n,p,q}$
with $M_{r,r_1,n,p,q}=\frac{[n]_{p,q}}{p^n}~\|L_{n,p,q}(f)-f\|_r>0$ $(\mbox{if}~\|L_{n,p,q}(f)-f\|_r)~\mbox{would be equal to}~0,$
this would imply that $f$ is a linear function, a contradiction).

Therefore, finally we get $\|L_{n,p,q}(f)-f\|_r~\geq~\frac{p^n}{[n]_{p,q}}~C_{r,r_1,p,q}(f)$ for all $n\in\mathbb{N},$ where
\begin{equation*}
  C_{r,r_1,p,q}(f)=\min\left\{M_{r,r_1,1,p,q}(f),...,M_{r,r_1,n_{\circ},p,q}(f),~\frac{1}{2}\|S_{p,q}(f)\|\right\},
\end{equation*}
which completes the proof.

Combining Theorem  \ref{theorem2.3}, and Theorem \ref{theorem2.1} (i), we immediately get the following result.

\begin{corollary}\label{corollary2.1} {\em Let $R>p^4/q^4>1,~f:D_R\rightarrow\mathbb{C}$
be analytic in $D_R$, i.e., $f(z)=\sum\limits_{k=0}^{\infty}c_kz^k,$ for all $z\in D_R$ and
let $1\leq r<\frac{q^3r_1}{p^3}<\frac{q^4R}{p^4}$ be arbitrary fixed. If $f$ is not a polynomial
of degree $\leq 1,$ then for all $n\in\mathbb{N}~\mbox{and}~|z|\leq r,$ we have
\begin{equation*}
   \|L_{n,p,q}(f)-f\|_r\sim\frac{p^n}{[n]_{p,q}},
\end{equation*}
where the constants in the equivalence depend  only on $f,~r,~r_1.~ p$ and $g$ but are independent of $n$.}
\end{corollary}

\section{Approximation results}

\hspace{8mm}Concerning the simultaneous approximation, we prove the following:
\begin{theorem}\label{theorem3.1}{\em Let $R>p^4/q^4>1,~f:D_R\rightarrow\mathbb{C}$
be analytic in $D_R$, i.e., $f(z)=\sum\limits_{k=0}^{\infty}c_kz^k,$ for all $z\in D_R$ and
let $1\leq r<r^{\ast}<\frac{p^3r_1}{q^3}<\frac{p^4R}{q^4}$ be arbitrary fixed. Also let $m\in\mathbb{N}.$
If $f$ is not a polynomial of degree $\leq\max\{1,m-1\},$ then for all $n\in\mathbb{N},$ we have
\begin{equation*}
   \|L_{n,p,q}^{(m)}(f)-f^{(m)}\|_r\sim\frac{p^n}{[n]_{p,q}},
\end{equation*}
where the constants in the equivalence depend  only on $f,~r,~r_1.~m,~p$ and $q$ but are independent of $n$.}
\end{theorem}

\noindent\textbf{Proof.} We already have the upper estimate for $\|L_{n,p,q}^{(m)}(f)-f^{(m)}\|_r,$
by Theorem \ref{theorem2.1} (ii), so it remains to find the lower estimate for $\|L_{n,p,q}^{(m)}(f)-f^{(m)}\|_r,$

Let us  denote by $\Gamma$ the circle of the radius $r^{\ast}$ and center $0.$
We have that the inequality $|\upsilon-z|\geq r^{\ast}-r$ holds for all $|z|\leq r$ and
$\upsilon\in\Gamma.$ Cauchy's formula is expressed by
\begin{equation}\label{e3.1}
\big| L_{n,p,q}^{(m)}(f)(z)-f^{(m)}(z)\big| =\frac{m!}{2\pi}\left|\int_{\gamma} \frac{L_{n,p,q}(f)(\upsilon)-f(\upsilon)}{(\upsilon-z)^{m+1}}d\upsilon\right|
\end{equation}
Now, as in the proof of Theorem \ref{theorem2.1} (ii), for all $\upsilon\in\Gamma$ and $n\in\mathbb{N},$
we have\\
$L_{n,p,q}(f)(z)-f(z)$
\begin{equation}\label{e3.2}
  \hspace{2cm}=\frac{p^n}{[n]_{p,q}}\left\{-S_{p,q}(f)(z)+\frac{p^n}{[n]_{p,q}}\left[\frac{[n]^2_{p,q}}{p^{2n}}
  \left(L_{n,p,q}(f)(z)-f(z)+\frac{p^n}{[n]_{p,q}}S_{p,q}(f)(z)\right)\right]\right\}
\end{equation}
By (\ref{e3.1}) and (\ref{e3.2}), we get
\begin{equation*}
  L_{n,p,q}^{(m)}(f)-f^{(m)}(f)=\frac{p^n}{[n]_{p,q}}\bigg\{~\frac{m!}{2\pi i}\int_{\Gamma}-~\frac{S_{p,q}(f)(z)}{(\upsilon-z)^{m+1}}d\upsilon\hspace{4cm}
\end{equation*}
\begin{equation*}
  ~+~\frac{p^n}{[n]_{p,q}}.\frac{m!}{2\pi i}\int_{\Gamma} ~\frac{[n]_{p,q}^2\left(L_{n,p,q}(f)(z)-f(z)+p^n\frac{S_{p,q}(f)(z)}{[n]_{p,q}}\right)}{p^{2n}(\upsilon-z)^{m+1}}d\upsilon\bigg\}
\end{equation*}
\begin{equation*}
=\frac{p^n}{[n]_{p,q}}\bigg\{[-S_{p,q}(f)(z)]{^{(m)}}~+~\frac{p^n}{[n]_{p,q}}.\frac{m!}{2\pi i}
  ~\int_{\Gamma}~\frac{[n]_{p,q}^2\left(L_{n,p,q}(f)(z)-f(z)+p^n\frac{S_{p,q}(f)(z)}{[n]_{p,q}}\right)}{p^{2n}(\upsilon-z)^{m+1}}d\upsilon\bigg\}.
\end{equation*}
Hence
\begin{equation*}
  \|L_{n,p,q}^{(m)}-f^{(m)}\|_r\geq\frac{p^n}{[n]_{p,q}}\bigg\{\left\|-\big[S_{p,q}(f)\big]^{(m)}\right\|_r\hspace{4cm}
\end{equation*}
\begin{equation*}
  ~-~\frac{p^n}{[n]_{p,q}}\bigg\|\frac{m!}{2\pi}\int_{\Gamma} ~\frac{[n]_{p,q}^2\left(L_{n,p,q}(f)(z)-f(z)+p^n\frac{S_{p,q}(f)(z)}{[n]_{p,q}}\right)}{p^{2n}(\upsilon-z)^{m+1}}d\upsilon\bigg\|_r\bigg\}.
\end{equation*}
Now by using Theorem \ref{theorem2.2} for all $n\in\mathbb{N},$ we get
\begin{equation*}
  \bigg\|\frac{m!}{2\pi}\int_{\Gamma}~\frac{[n]_{p,q}^2\left(L_{n,p,q}(f)(z)-f(z)
  +p^n\frac{S_{p,q}(f)(z)}{[n]_{p,q}}\right)}{p^{2n}(\upsilon-z)^{m+1}}d\upsilon\bigg\|_r\hspace{3cm}
\end{equation*}
\begin{equation*}
  \leq\frac{m!}{2\pi}\frac{2\pi r^{\ast}[n]_{p,q}^2}{(r^{\ast}-r)^{m+1}p^{2n}}\left\|L_{n,p,q}(f)-f
  +p^n\frac{S_{p,q}(f)}{[n]_{p,q}}\right\|_{r^{\ast}}
\end{equation*}
\begin{equation*}
  \leq Q_{r_1,p,q}(f).\frac{m!~r^{\ast}}{(r^{\ast}-r)^{m+1}}.\hspace{5.6cm}
\end{equation*}
\hspace{8mm}But by hypothesis on $f,$ we have $\big\|-\big[S_{p,q}(f)\big]^{(m)}\big\|_{r^{\ast}}>0.$
Indeed, supposing the contrary, it would follow that $\big[S_{p,q}(f)\big]^{(m)}(z)=0,$ for all
$|z|\leq r^{\ast},$ where by the statement of Theorem \ref{theorem2.2}, we have
\begin{equation*}
S_{p,q}(f)(z)=\sum\limits_{k=2}^{\infty}q^{n-(k-1)}c_k\frac{[k]_{p,q}-[k]_q}{q-1}z^k
=\sum\limits_{k=2}^{\infty}q^{n-(k-1)}c_k\big([1]_{p,q}+...+[k-1]_{p,q}\big)z^k.
\end{equation*}
Firstly, supposing that $m=1,$ by $S^{\prime}_{p,q}(f)(z)=\sum\limits_{k=2}^{\infty}q^{n-(k-1)}c_k~k~\big([1]_{p,q}+...+[k-1]_{p,q}\big)z^{k-1}=0,$
for all $|z|\leq r^{\ast},$ would follow that $c_k=0,$ for all $k\geq 2,$ that is, $f$ would
be a polynomial of degree $1=\max\{1,m-1\},$ a contradiction with the hypothesis.

Taking $m=2,$ we would get $S^{\prime\prime}_{p,q}(f)(z)=\sum\limits_{k=2}^{\infty}q^{n-(k-1)}~c_k~k~(k-1)
~\big([1]_{p,q}+...+[k-1]_{p,q}\big)z^{k-2}=0,$ for all $|z|\leq r^{\ast},$ which immediately would imply that $c_k=0,$
for all $k\geq 2,$ that is, $f$ would be a polynomial of degree $1=\max\{1,m-1\},$ a contradiction with the hypothesis.

Now, taking $m>2,$ for all $|z|\leq r^{\ast},$ we would get
\begin{equation*}
S^{(m)}_{p,q}(f)(z)=\sum\limits_{k=m}^{\infty}q^{n-(k-1)}~c_k~k~(k-1)...(k-m+1)
~\big([1]_{p,q}+...+[k-1]_{p,q}\big)z^{k-m}=0,
\end{equation*}
which would imply that $c_k=0,$ for all $k\geq m,$ that is, $f$ would be a polynomial of degree
$m-1=\max\{1,m-1\},$ a contradiction with the hypothesis.

Finally, we prove some approximation results for the iterates of $(p,q)$-Lorentz operators.

For $f$ analytic in $D_R$ that is of the form $f(z)=\sum\limits_{k=0}^{\infty}c_kz^k,$ for all $z\in D_R,$
let us define the iterates of complex Lorentz operators $L_{n,p,q}(f)(z),$ by $L_{n,p,q}^{(1)}(f)(z)=L_{n,p,q}(f)(z)$
and $L_{n,p,q}^{(m)}(f)(z)=L_{n,p,q}[L_{n,p,q}^{(m-1)}(f)](z),$ for any $m\in\mathbb{N},~m\geq 2.$

Since we have $L_{n,p,q}(f)(z)=\sum\limits_{k=0}^{\infty}c_kL_{n,p,q}(e_k)(z),$ by
recurrence for all $m\geq 1,$ we get that $L_{n,p,q}^{(m)}(f)(z)=
\sum\limits_{k=0}^{\infty}c_kL^{(m)}_{n,p,q}(e_k)(z),$
where $L_{n,p,q}^{(m)}(e_k)(z)=1,~\mbox{if}~k=0,~L_{n,p,q}^{(m)}(e_k)(z)=z
~\mbox{if}~k=1,~L_{n,p,q}^{(m)}(e_k)(z)=0,~\mbox{if}~~k\geq n+1$ and
 \begin{equation*}
 L_{n,p,q}^{(m)}(e_j)(z)=\left(1-p^{n-1}\frac{[1]_{p,q}}{[n]_{p,q}}\right)^m
\left(1-p^{n-2}\frac{[2]_{p,q}}{[n]_{p,q}}\right)^m...\left(1-p^{n-(j-1)}\frac{[j-1]_{p,q}}{[n]_{p,q}}\right)^m~z^k,
\end{equation*}
for $2\leq k\leq n.$

We present the following:
\begin{theorem}\label{theorem3.2} {\em Let $R>p>q>1$ and $1\leq r < \frac{pr_1}{q}<\frac{pR}{q}$ be arbitrary fixed.
Denoting $D_{R}=\{z\in\mathbb{C}: |z|<R\}.$ Suppose that $f:D_R\rightarrow\mathbb{C}$ is analytic in $D_R$, i.e.,
$f(z)=\sum\limits_{k=0}^{\infty}c_kz^k,$ for all $z\in D_R,$  we have the upper estimate
 \begin{equation*}
\left\| L_{n,p,q}^{(m)}(f)-f\right\|_r \leq \frac{mp^n}{[n]_{p,q}}\frac{q-p+1}{(q-p)^2}\sum\limits_{k=0}^{\infty} |c_k| (k+1) r_1^k.
\end{equation*}
Moreover, if $\lim\limits_{n\rightarrow\infty}\frac{m_np^n}{[n]_{p,q}}=0,$ then
 \begin{equation*}
\lim\limits_{n\rightarrow\infty}\left\| L_{n,p,q}^{(m_n)}(f)-f\right\|_r =0.
\end{equation*}}
\end{theorem}
\noindent\textbf{Proof.} For all $|z|\leq r,$ we easily obtain

\vspace{3mm}
$\big|f(z)-L^{(m)}_{n,p,q}(f)(z)\big|$
\begin{equation*}
  \leq\sum\limits_{k=2}^n|c_k|r^k\left[1-\left(1-p^{n-1}\frac{[1]_{p,q}}{[n]_{p,q}}\right)^m
...\left(1-p^{n-(j-1)}\frac{[j-1]_{p,q}}{[n]_{p,q}}\right)^m\right]+\sum\limits_{k=n+1}^{\infty}|c_k|r^k.
\end{equation*}
Denoting $A_{k,n}=\left(1-p^{n-1}\frac{[1]_{p,q}}{[n]_{p,q}}\right)
...\left(1-p^{n-(j-1)}\frac{[k-1]_{p,q}}{[n]_{p,q}}\right),$ we get $1-A^m_{k,n}=\\
(1-A_{k,n})(1+A_{k,n}+A_{k,n}^2+...+A_{k,n}^{m-1})\leq m(1-A_{k,n})$ and therefore since
$1-A_{k,n}\leq p^{n-(k-1)}\frac{(k-1)[k-1]_{p,q}}{[n]_{p,q}},$ for all $|z|\leq r,$ we obtain
\begin{equation*}
  \sum\limits_{k=2}^n|c_k|r^k\left[1-\left(1-p^{n-1}\frac{[1]_{p,q}}{[n]_{p,q}}\right)^m
...\left(1-p^{n-(j-1)}\frac{[j-1]_{p,q}}{[n]_{p,q}}\right)^m\right]+\sum\limits_{k=n+1}^{\infty}|c_k|r^k.
\end{equation*}
\begin{equation*}
  \leq m\sum\limits_{k=2}^{\infty}|c_k|r^kp^{1-(k-1)}~[1-A_{k,n}]\leq\frac{mp^{n+1}}{[n]_{p,q}}\sum\limits_{k=2}^{\infty}|c_k|r^k(k-1)[k-1]_{p,q}r^k
\end{equation*}
\begin{equation*}
  \leq\frac{mp^{n+1}}{[n]_{p,q}}\sum\limits_{k=2}^{\infty}|c_k|r^k\frac{kq^k/p^k}{q-p}\leq\frac{mp^{n+1}}{[n]_{p,q}}\frac{1}{q-p}
  \sum\limits_{k=2}^{\infty}|c_k|(k+1)(\frac{q}{p}r)^k\hspace{0.4cm}
\end{equation*}
\begin{equation*}
 \leq\frac{mp^{n+1}}{[n]_{p,q}}\frac{1}{q-p}
  \sum\limits_{k=2}^{\infty}|c_k|(k+1)(r_1)^k.\hspace{4.4cm}
\end{equation*}
On the other hand, following exactly the reasonings in the proof of the Theorem \ref{theorem2.1},
we get the estimate
\begin{equation*}
  \sum\limits_{k=n+1}^{\infty}|c_k|r^k\leq \frac{p^{n+1}}{[n]_{p,q}}\frac{\sum\limits_{k=0}^{\infty}|c_k|(k+1)(r_1)^k}{(q-p)^2}\leq \frac{mp^{n+1}}{[n]_{p,q}}.\frac{\sum\limits_{k=0}^{\infty}|c_k|(k+1)(r_1)^k}{(q-p)^2}.
\end{equation*}
Collecting now all the estimates and taking into account that $\frac{1}{q-p}+\frac{}{(q-p)^2}=\frac{q-p+1}{(q-p)^2},$
we arrive at the desired estimate.

Since $\lim\limits_{n\rightarrow\infty}\frac{m_np^n}{[n]_{p,q}}=0,$ it follows the conclusion that
 \begin{equation*}
\lim\limits_{n\rightarrow\infty}\left\| L_{n,p,q}^{(m_n)}(f)-f\right\|_r =0.
\end{equation*}

\end{document}